\documentclass[12pt]{amsart}
\usepackage{graphics}

\def\gap{\smallskip\noindent}
\def\mgap{\medskip\noindent}
\def\hils{{\mathcal H}}
\def\Re{{\rm Re}\,}
\def\Im{{\rm Im}\,}
\def\R{{\Bbb R}}
\def\Z{{\Bbb Z}}

\def\C{{\Bbb C}}

\begin{document}

\title[Solutions, spectrum, and dynamics for Schr\"odinger operators]
{Solutions, spectrum, and dynamics for Schr\"odinger operators on infinite
domains}

\author{A.~Kiselev and Y.~Last}


\begin{abstract}

Let $H^{\Omega}_{V}= -\Delta +  V(x)$ be a Schr\"odinger operator defined on
an unbounded domain $\Omega\subset\R^d$ with Dirichlet boundary
conditions on $\partial \Omega$ (when $\Omega=\R^{d}$ there is
no boundary condition, of course). Let $u(x,E)$ be a solution of the
Schr\"odinger equation
$(H^{\Omega}_{V}-E)u(x,E)=0$, and let $B_{R}$ denote a ball of radius $R$
centered at zero.
We show relations between the rate of growth of the $L^{2}$
norm $\|u(x,E)\|_{L^{2}(B_{R}\cap \Omega)}$ of such solutions,
as $R \to \infty$, and continuity properties of spectral measures of the
operator $H_{V}^{\Omega}$.
These results naturally lead to new criteria for the identification
of various spectral properties.
We also prove new fundamental relations between the rate of growth of $L^{2}$
norms of generalized eigenfunctions, dimensional properties of spectral
measures, and dynamical properties of the corresponding quantum systems.
We apply these results to study transport properties of some particular
Schr\"odinger operators.

\end{abstract}

\maketitle

\begin{center}
\bf \large 1. Introduction and main results
\end{center}
\mgap
In this paper we investigate the relations between the rate of decay of
solutions of Schr\"odinger equations, continuity properties of
spectral measures of the corresponding operators, and dynamical properties
of the corresponding quantum systems.
The first main result of this paper shows that, in great
generality, certain upper bounds on the rate of growth of
$L^{2}$ norms of generalized
eigenfunctions over expanding balls imply certain minimal singularity
of the spectral measures.  Consider an
operator $H_{V}^{\Omega}$ defined by the differential expression
\[  H^{\Omega}_{V} = -\Delta + V(x) \]
on some connected infinite domain $\Omega$ with a smooth boundary
and with Dirichlet boundary conditions on $\partial \Omega$.
The case of $\Omega = \R^{d}$ is not excluded; no boundary conditions
are needed in this case.
To every vector $\phi \in L^{2}(\Omega)$ we associate a spectral measure
$\mu^{\phi}$ in the usual way
(namely, $\mu^{\phi}$ is the unique Borel measure on $\R$ obeying
$\int f(E)\,d\mu^{\phi}(E) = (f(H^{\Omega}_{V})\phi, \phi)$ for any
Borel function $f$).
For any measure $\mu$, we define the upper $\alpha$-derivative
$D^{\alpha}\mu(E)$ in the standard way:
\[ D^{\alpha}\mu(E) = \limsup_{\delta \rightarrow 0}
\frac{\mu(E-\delta, E+\delta)}{\delta^{\alpha}}. \]
We denote by $B_{R}$ the ball of radius $R$ centered at the origin, and
use the notation $\|f\|_{B_{R}}$ for the $L^{2}$ norm of the function $f$
restricted to $B_{R}$.  We  
 denote by $W_{m}^{l}$ the usual Sobolev
spaces of functions $f$ such that $D^{l}f$ exists in the distributional sense 
and $\int(|u|^{m}+|D^{l}u|^{m})dx < \infty.$ We say that $f(x) \in W^{l}_{m, {\rm loc}}(\Omega)$
if $f(x) \in W^{l}_{m}(\Omega \cap B_R)$ for every $R <\infty.$ 
One of the main theorems that we prove here is the following:

\bigskip\noindent
\bf Theorem 1.1. \it Assume that the potential $V(x)$ belongs to $L^{\infty}_{{\rm loc}}$ 
and is bounded from below, and $\Omega$ is a domain with piecewise smooth boundary.
  Suppose that there exists a distributional solution $u(x, E)$ of the generalized
eigenfunction equation
\begin{equation}
\label{1.1}
 (H^{\Omega}_{V}-E)u(x,E)=0
\end{equation}
satisfying the boundary conditions and
such that for
some $\alpha,$ $0 \leq \alpha \leq 1$, we have
\begin{equation}
\label{1.2}
\liminf_{R \to \infty} R^{-\alpha}\int\limits_{B_{R}\cap \Omega}|
u(x,E)|^{2}\,dx < \infty.
\end{equation}
Fix some compactly supported $\phi(x) \in L^{2}(\Omega)$ such that
\[ \int\limits_{\Omega}\phi(x) u(x,E)\,dx \neq 0. \]
Then we have
\[ D^{\alpha}\mu^{\phi}(E) >0. \]

\bigskip\noindent
Remarks. \rm 1. Notice that under our assumptions on the potential, we have 
$u \in W^{2}_{2,{\rm loc}}$
by standard results on Sobolev estimates for elliptic operators (see, e.g., \cite{Ev}),
and the boundary values for $u$ are well-defined.

\smallskip
2. We chose not to formulate this Theorem for more general classes of potentials,
domains, and boundary conditions in order to be able to give a transparent proof.
Certainly, we can extend this theorem to wider classes
of potentials and boundary conditions. The nature of the limitations
will be clear from the proof and Stark operators example in Appendix 2.
For instance, when $\Omega = \R^{d}$, we only ask that the negative part of
the potential, $V_{-}$, belongs to the Kato class $K^{d}$ 
(see, e.g.,
\cite{AS,Sim1} for the definition of Kato classes).

\smallskip
3. If we replace ``$<\infty$'' in \eqref{1.2} by ``$= 0$'', we
obtain that $D^{\alpha}\mu^{\phi}(E) = \infty$.\\

 \rm  Theorem 1.1 provides information on the pointwise behavior
of spectral measures from rather simple and natural assumptions
about the behavior of generalized eigenfunctions. From this theorem
follow new criteria for the existence of absolutely continuous
spectrum or singular continuous spectrum of given dimensional characteristics
(see Section 2, and, in particular, Theorem 2.5 for more details).
This contrasts the well-known result
(see \cite{Ber,Schn,Sim1}) that existence of a
polynomially bounded (but not $L^2$)
solution of (1) implies that the energy $E$ belongs to the essential
spectrum of $H_V,$ but gives no further information on the structure
of the essential spectrum.
To the best of our knowledge, Theorem 1.1 is the first rigorous result
providing a relation between the behavior of solutions and pointwise
properties of the spectral measures for multidimensional Schr\"odinger
operators.

A result analogous to Theorem 1.1 also holds for discrete Schr\"odinger
operators defined on some $\Omega \subset \Z^{d}$ by
\[ (h_{v}u)(n) = \sum\limits_{|m-n|=1,\, m \in \Omega} u(m) +v(n)u(n). \]
We discuss this extension in Section 3.
In Appendix 1, we also indicate that  results similar
to Theorem 1.1 hold for more general elliptic
and higher order operators.

The motivation for seeking relations between the pointwise in energy
behavior of solutions and properties of spectral measures comes from the
fact that in many problems the solutions are among the objects we can
hope to investigate. When we are
interested in the fine structure of the spectrum of Schr\"odinger
operators for which the methods of scattering theory are not applicable,
there are very limited tools in higher dimensions which may be
effectively used for spectral analysis. On the other hand, for
one-dimensional Schr\"odinger operators the subordinacy theory created
by Gilbert and Pearson \cite{GiPe,Gi} and further extended by
Jitomirskaya and Last \cite{JiLa,JiLai,JiLaii}
provides a powerful method for spectral analysis.
The main results of the above mentioned papers give a necessary and
sufficient link between the behavior of solutions and the singularity
of the spectral measure.
Subordinacy theory played an important role in many recent results in
one-dimensional spectral theory
\cite{CK,Da,Ji,JiLa,JiLai,JiLaii,Ki,KLS,LaSi,Re}.

In this paper, we derive only a sufficient-type relation between the
solutions and the spectrum, but in much greater generality.
However, in contrast to subordinacy theory, which requires
comparison of different solutions, we need information about only one
solution---the one obeying the appropriate boundary conditions.
We remark that for one-dimensional Schr\"odinger operators, the result of
Theorem 1.1 can be derived from subordinacy theory \cite{JiLai,JiLaii}.

Our second major result in this paper establishes a fundamental relation
between spectral properties, generalized eigenfunctions and quantum dynamics,
and in particular, provides new bounds for the
transport properties of quantum systems. We study the behavior of the
time-averaged moments of the position operator $X$ under the
Schr\"odinger evolution.
Pick some initial state $\psi$ and consider
\[ \langle \langle |X|^{m}\rangle \rangle_{T} =
\frac{1}{T} \int\limits_{0}^{T}|\langle |X|^{m} \exp (-iH^{\Omega}_{V}t)\psi,
\exp (-iH^{\Omega}_{V}t) \psi \rangle |\,dt. \]
Recall that a measure $\mu$ is called $\alpha$-continuous if
it gives zero weight to any set of zero $\alpha$-dimensional
Hausdorff measure
(we recall the definition of these measures in Section 2).
Let us  denote by $P_{\alpha c}$ the spectral projector on the
$\alpha-$continuous spectral subspace, the set of all vectors
$\xi$ such that $\mu^{\xi}$ is $\alpha$-continuous
(see \cite{L}).
In particular, if $\mu^{\psi}$ has an $\alpha$-continuous
component (i.e., $P_{\alpha c}\psi \ne 0$), then the
following lower bound holds \cite{Comb,Gua,Gub,L}:
\[ \langle \langle |X|^{m}\rangle \rangle_{T} \geq C_m T^{\frac{m\alpha}{d}} \]
(here $d$ is the space dimension and $C_m$ is a constant depending on 
$\mu^{\psi}$ and $m$).

Recall that for a wide class of Schr\"odinger operators, one has
a generalized eigenfunction expansion theorem
(see, e.g., \cite{Ber,LaSi,Sim1}).
In particular, for every $\psi$ there is a unique unitary map
$U_{\psi}$ from the
cyclic subspace $\hils_{\psi}$, generated by the vector $\psi$ and
the operator $H_{V}$,
to  $L^{2}(\R, d \mu^{\psi}(E))$. This map sends $\psi$ to a function
equal to $1$ everywhere and realizes a unitary equivalence
$U_{\psi}H_{V}|_{\hils_{\psi}}U^{-1}_{\psi}=E,$
where $E$ stands for the operator of multiplication by $E.$ The operator
$U_{\psi}$ is an integral operator with kernel $u(x, E),$ where the $u(x,E)$'s,
for each fixed $E$, solve \eqref{1.1} and are called generalized
eigenfunctions. 
We will say that the $u(x,E)$'s correspond to $\psi$ if they constitute
the kernel of the unitary map $U_{\psi}$ described above. Note that they
are only defined a.e.\ w.r.t. $\mu^\psi$.
We prove the following theorem, which holds in both the discrete
and continuous settings:

\bigskip\noindent
\bf Theorem 1.2. \it Let $\psi$ be a vector for which there exists
a Borel set $S\subset\R$ of positive $\mu^\psi$ measure, such that the
restriction of $\mu^\psi$ to $S$ is $\alpha$-continuous
and, in addition, the generalized eigenfunctions $u(x,E)$
for all $E\in S$ satisfy
\begin{equation}
\label{1.3}
\limsup_{R \to \infty} R^{-\gamma}\|u(x,E)\|^{2}_{B_{R}} < \infty
\end{equation}
for some $\gamma$ such that $0 < \gamma < d$.
Then, for any $m > 0$, there exists a constant $C_m$ such that
\begin{equation}
\label{1.4}
\langle \langle |X|^{m}\rangle \rangle_{T} \geq C_m T^{\frac{m \alpha}{\gamma}}
\end{equation}
for all $T>0$.

\bigskip\noindent
Remarks. \rm
1. Theorem 1.2 is somewhat related to (although it does not coincide
with) some recent heuristic results by Ketzmerick et. al. \cite{Ketz}.

\smallskip
2. It may be seen from Theorem 1.1 that we cannot have $\gamma < \alpha,$ since it
would follow that the upper $\gamma$-derivative of the spectral measure is
positive on too large a set
(see Corollary 2.6). The physical reason is that when $V$ is bounded from 
below, the velocity is bounded, and the propagation rate is at most ballistic. 
However, the range of applicability
of Theorem 1.2 is wider than that of Theorem 1.1.
In particular, it is applicable
to operators with strongly negative potentials, such as Stark operators, which
exhibit faster-than-ballistic transport. See Appendix 2. \\

The somewhat striking aspect
of Theorem 1.2 is that for a fixed non-zero
spectral dimension, faster decay of $u(x,E)$
leads to faster transport. Theorem 1.2 shows that the
behavior of the generalized eigenfunctions plays an important role in
determining dynamical properties of quantum systems.
We apply Theorem 1.2 to investigate the dynamics in the
random decaying potentials model studied in \cite{KLS}.
When weakly coupled, these systems have (almost surely) some singular
continuous spectrum with local dimensions that depend on the energy,
but we show that the
dynamical spreading of wavepackets, for any energy region where the
spectrum is continuous, is almost ballistic with probability one.
More precisely, we show that for almost
every realization, we have for every $\epsilon > 0$
a bound of the form
\[   \langle \langle |X|^{m}\rangle \rangle_{T} \geq
C_{m,\epsilon}T^{m(1-\epsilon)}. \]

The paper is organized as follows. In Section 2 we prove
Theorem 1.1
 and its corollaries, rendering new spectral criteria.
 In Section 3 we sketch the
argument for similar results in the discrete setting.
In Section 4 we consider some
simple examples, in particular, showing that the result of Theorem 1.1
provides only a sufficient but not necessary criterion for positivity of the
derivative of the spectral measure. It is, however, an optimal result in the
sense that one cannot in general say more by looking only at the
rate of growth of the $L^{2}$ norm (Section 5).
 It remains an interesting open
question to find additional properties of solutions
that determine the spectrum (or other important characteristics of the
operator, such as transport properties) completely.
In Section 5 we study the relationship between solutions, spectral
dimension, and quantum dynamics, in particular proving Theorem 1.2.
 In the appendices, we indicate further possible
generalizations for elliptic and higher order operators and
consider dynamics for strongly perturbed one-dimensional Stark operators.
The example of Stark operators provides another illustration of the
relationship between the behavior of solutions and transport properties.
\gap
\begin{center}
\large \bf 2. Solutions and spectrum: Continuous case
\end{center}
\mgap
\rm We begin the proof of Theorem 1.1 with the following
simple observation: \\

\gap
\bf Lemma 2.1 \it Let $A$ be a self-adjoint operator acting on a Hilbert
space ${\mathcal H}$ and fix a vector $\phi \in {\mathcal H}$.
Let $z \in \C \setminus \R$. Then
\[  \Im z \| (A-z)^{-1}\phi \|^{2}_{\mathcal H} = \Im ((A-z)^{-1}\phi,
 \phi). \]
\gap
\begin{proof} Consider the spectral representation associated
with a vector $\phi$ and perform a straightforward computation:
\[  \Im\left(\int\limits_{R}\frac{d\mu^{\phi}}{t-z}\right)=
 \Im z \int \limits_{R}\frac{d\mu^{\phi}}{|t-z|^{2}}. \]
\end{proof}   \rm

The first idea in the proof of Theorem 1.1 is to estimate from below
$\Im((H_{V}^{\Omega}-E-i\epsilon)^{-1}\phi,\phi)$ as $\epsilon
\rightarrow 0.$
Such an estimate is equivalent to an estimate on the upper $\alpha$-derivative
of the
spectral measure by the following lemma: \\

\gap
\bf Lemma 2.2. \it Let $Q^{\beta}_{\mu}(E)$ denote
\[ Q^{\beta}_{\mu}(E)= \limsup_{\epsilon \rightarrow 0}
\epsilon^{\beta} \, \Im \left( \int \frac{d\mu(t)}{t-E-i\epsilon} \right). \]
Then
\[ D^{\alpha}\mu(E) \leq C_{1} Q^{1-\alpha}_{\mu}(E)
\leq C_{2}D^{\alpha}\mu(E), \]
where $C_{1},$ $C_{2}$ are positive constants depending only on $\alpha.$
\gap
\begin{proof} \rm The proof is a direct computation. For
details, we refer to \cite{DJLS}, Lemmas 3.2 and 3.3.
\end{proof}     \rm

\rm To derive an estimate on the imaginary part of the Borel transform, we
will use Lemma 2.1, namely  estimates from below on the norm of the
function
\[ \theta (x, E+i\epsilon) = (H_{V}^{\Omega}-E-i\epsilon)^{-1}\phi (x) \]
over balls of radius of order $\frac{1}{\epsilon}$ as $\epsilon$ goes
to zero over some properly chosen sequence.

The last technical lemmas that we need for the proof concern
estimation of the $W_{2}^{1}$ norms of $u(x,E)$ and $\theta(x, z)$
in terms of their $L^{2}$ norms. \\

\gap
\bf Lemma 2.3. \it Let $\Omega \subset \R^{d}$ be a domain with
piecewise smooth boundary.
Suppose that the potential $V$ belongs to $L^{\infty}_{{\rm loc}}$ and 
is bounded from below,
 and let $H_{V}^{\Omega}$ denote an operator with Dirichlet
boundary conditions on $\partial \Omega.$ Suppose that the function
$g(x,z)$ satisfies Dirichlet boundary conditions and
\[ (H_{V}^{\Omega}-z) g(x, z)= \phi(x), \]
where $\phi \in L^{2}(\Omega)$ is compactly supported and real-valued, and $z$ is
in general complex. Then
\begin{equation}
\label{2.1}
\|g\|_{W_{2}^{1}(B_{R} \cap \Omega)} \leq C(z, V_{-}) \left( \|g\|_{L^{2}
(B_{R+1} \cap \Omega)}+ \|\phi\|_{L^{2}(\Omega)} \right).
\end{equation}
The constant in \eqref{2.1}
 depends only on the lower bound on $V$ and on $z,$
and may be chosen uniformly for $z$ in any compact set.
\gap
\begin{proof} The proof is standard and we provide it for the sake
of completeness. See, for example, \cite{AS,Sim1} for detailed exposition 
of similar results and further references.
Throughout the proof, we assume
that the function $g$ is sufficiently smooth to justify integration by parts
(local $W_{2}^{2}$ is sufficient). Clearly this is the case under our assumptions
on $V$ (see, e.g., \cite{Ev}).
To prove the bound \eqref{2.1} with the constant independent
of $R,$ let
\[ g(x,z) = g_{1}(x, z) + i g_{2}(x, z), \]
where $g_{1},$ $g_{2}$ are real-valued. For any $\psi \in C^{\infty}(\Omega)$
such that $1 \geq \psi(x) \geq 0,$ $\psi(x) =1$  when $x \in B_{R} \cap
\Omega,$ $\psi(x) =0$ when $x \notin B_{R+1} \cap \Omega,$ we have
\begin{eqnarray}
\int\limits_{B_{R} \cap \Omega} (\nabla g_{1})^{2}\,dx & \leq &
\int\limits_{B_{R+1} \cap \Omega} \psi (\nabla g_{1})^{2}\,dx \,=\,
\int\limits_{\partial (\Omega \cap B_{R+1})} \psi \frac{\partial g_{1}}{\partial
n} g_{1}\, d\sigma - \nonumber \\
 & &- \int\limits_{B_{R+1} \cap \Omega} (\nabla \psi)
(\nabla g_{1}) g_{1}\,dx
 -  \int\limits_{B_{R+1} \cap \Omega}
\psi g_{1}\Delta g_{1}\,dx,
\label{2.2}
\end{eqnarray}
where $d \sigma$ is the surface measure on $\partial (\Omega \cap B_{R+1})$
induced from $\R^{d}.$  The first term vanishes because $g_{1}$
vanishes on $\partial \Omega$ and $\psi$ vanishes on $(\partial B_{R}) \cap
\Omega.$ Furthermore, by Green's formula
\begin{equation}
\label{2.3}
 2 \int\limits_{B_{R+1} \cap \Omega} (\nabla \psi) (\nabla g_{1}) g_{1}\,dx=
 \int\limits_{\partial(B_{R+1} \cap \Omega)} \frac{\partial\psi}{\partial n}
 (g_{1})^{2}\,d \sigma -
\int\limits_{B_{R+1} \cap \Omega} \Delta\psi ( g_{1})^{2}\,dx.
\end{equation}
The boundary term in this equality  is also equal to zero.
Substituting \eqref{2.3} into \eqref{2.2}, we find
\begin{eqnarray*} \int\limits_{B_{R} \cap \Omega}  (\nabla g_{1})^{2}\,dx
& \leq &
 \frac{1}{2}
\int\limits_{B_{R+1} \cap \Omega} \Delta \psi ( g_{1})^{2}\,dx \\
& & + \int\limits_{B_{R+1} \cap \Omega} \psi g_{1}\left((\Re z-V) g_{1}
+\phi - (\Im z) g_{2}\right)\,dx.
\end{eqnarray*}
Therefore,
\[ \|g_{1}\|^{2}_{W_{2}^{1}(B_{R})} \leq C_{\psi} \left( \|\phi\|^{2}_{L^{2}}
+(2(1+|z|) +\|V_{-}\|_{L^{\infty}})\|g_{1}\|^{2}_{L^{2}(B_{R+1})}+
 \Im z \|g_{2}\|^{2}_{L^{2}(B_{R+1})} \right). \]
A similar estimate holds for $g_{2}.$ Combining these two estimates, we obtain
the result of the lemma.
\end{proof} \rm

\noindent \it Remarks. \rm 1. We have not tried to determine the most
general classes
of potentials and boundary conditions for which the lemma holds.
With slightly more technical effort, we can treat
some other boundary conditions, such as Neumann, for instance.  

\smallskip
 2. For the case of the whole space, the lemma is true under the
assumption that $V_{-} \in K^{d},$ the Kato class, which allows singularities
in the negative part of the potential
(see \cite{Sim1} for the definition and properties
of potentials from these classes). This result follows from the technique
developed in \cite{AS, Sim1}, which uses Brownian motion to derive
subsolution estimates implying bounds like in Lemma 2.3. Although
\cite{AS,Sim1} consider only real $z$
(and homogeneous equation),
it is not hard to see that their arguments extend to give results like
\eqref{2.1}.

\smallskip
We now introduce an important object in our consideration.
Suppose $S$ is a domain with piecewise 
smooth boundary and $f,$ $g$ belong to
$W_{2,{\rm loc}}^{2}(S).$   We
denote by
$W_{\partial S}[f,g]$ the following expression
\begin{equation}
\label{2.4}
 W_{\partial S}[f,g] = \int\limits_{\partial S} (f(t)\frac{\partial g}
{\partial n}(t) - \frac{\partial f}{\partial n}(t)g(t))d\sigma(t),
\end{equation}
where $\sigma$ is the surface measure induced from $\R^{d}$ and
$\frac{\partial}{\partial n}$ is the derivative in the outer 
normal direction.
The definition makes sense for $W^{2}_{2, {\rm loc}}$ functions
by Sobolev trace theorems (see, e.g., \cite{GT}).
The notation $W$
stresses the fact that in one dimension, the corresponding expression is
related to the Wronskian of two functions (precisely, it is 
the difference of the
Wronskians taken at the endpoints of the interval $S$).
We will abuse verbal notation and call the expression \eqref{2.4}
the Wronskian of $f$ and $g$ over $\partial S$ for the rest of this paper. 
 The final lemma we need is \\

\gap
\bf Lemma 2.4. \it Suppose that two functions $f,$ $g$ are
locally $W_{2}^{2}$ and satisfy Dirichlet boundary condition
on $\partial \Omega.$
Then for every $R$
\[ \int\limits_{0}^{R} |W_{\partial (B_{r} \cap \Omega)}[f,g]| \, dr \leq
\|f\|_{W_{2}^{1}(B_{R}\cap \Omega)}\|g\|_{W_{2}^{1} (B_{R}\cap \Omega)}. \]
\gap
\begin{proof} \rm We have $W_{\partial \Omega \cap B_{R}}[f,g]=0$
since $f$ and $g$ satisfy the boundary conditions.
Next note that
\[  \int\limits_{0}^{R} |W_{\partial B_{r} \cap \Omega )}[f,g]| \, dr \leq
\int\limits_{B_{R}\cap \Omega}(|f||\nabla g|+|\nabla f||g|)\,dx \leq
\|f\|_{W_{2}^{1}(B_{R}\cap \Omega)}\|g\|_{W_{2}^{1} (B_{R} \cap \Omega)}. \]
We used the Cauchy-Schwartz inequality in the last step.
\end{proof}

\rm Now we are ready to prove Theorem 1.1.

\begin{proof} \rm
 An interplay of the scales in space and in the spectral parameter plays an
important role in the analysis. Let us assume that
\[ \int\limits_{\Omega} \phi(x)u(x, E)\,dx = c \ne 0. \]
Take sufficiently large $R_{0},$ such that ${\rm supp} \phi \subset
B_{R_{0}}.$
By Green's formula we have
\[  E \int\limits_{B_{R_{0}}\cap \Omega} \theta (x, E+i\epsilon)
u(x,E)\,dx =
W_{\partial (B_{R_{0}}\cap \Omega)}[\theta, u] +
\int\limits_{B_{R_{0}}\cap \Omega} H_{V}^{\Omega}\theta (x, E+i\epsilon)
u(x,E)\,dx = \]
\[ = W_{\partial (B_{R_{0}}\cap \Omega)}[\theta, u] + (E+i\epsilon)\int
\limits_{B_{R_{0}}\cap \Omega} \theta(x, E+i\epsilon) u(x,E)\,dx+
\int\limits_{B_{R_{0}}\cap \Omega} \phi(x) u(x,E)\,dx. \]
 In the above computation we used the definition of
$\theta (x, z)$ and the fact that the function $u$ satisfies
$(H_{V}^{\Omega}-E)u=0.$ Hence, we obtain
\begin{equation}
\label{2.5}
W_{\partial (B_{R_{0}}\cap \Omega)}[\theta, u] = -c -i\epsilon
\int\limits_{B_{R_{0}}\cap \Omega} \theta(x, E+i\epsilon) u(x,E)\,dx.
\end{equation}
Let us integrate equation \eqref{2.5} from $R_{0}$ to some larger value
of $R:$
\[ \int\limits_{R_{0}}^{R} |W_{\partial (B_{r})\cap \Omega}[\theta, u]|\,dr
\geq |c|(R-R_{0}) - \epsilon \int\limits_{R_{0}}^{R} dr \left|\,
\int\limits_{B_{r}\cap \Omega} \theta(x, E+i\epsilon) u(x,E)\,dx \right|.
\]
Using Lemmas 2.3, 2.4, we see that
\begin{eqnarray}
\label{2.6}
C^{2}(\|\theta(x,E+i\epsilon)\|_{L_{2}(B_{R+1}\cap \Omega)}+\|\phi\|_{L^{2}})\|u\|_{L_{2}
(B_{R+1}\cap \Omega)} \geq \nonumber \\
 |c|(R-R_{0}) -\epsilon \int\limits_{0}^{R} dr
\|\theta(x,E+i\epsilon)\|_{L_{2}(B_{r}\cap \Omega)}
\|u\|_{L_{2}(B_{r}\cap \Omega)}. 
\end{eqnarray}
According to the assumption \eqref{1.2}
 of the theorem, there exists a sequence
$R_{n} \rightarrow \infty,$ such that
\begin{equation}
\label{2.7}
\|u\|_{L_{2}(B_{R_{n}}\cap \Omega)} \leq C_{1}R_{n}^{\frac{\alpha}{2}}.
\end{equation}
Let us set $\epsilon_{n}= \frac{C_{2}}{R_{n}},$ and pick
$R+1=R_{n}$ and $\epsilon = \epsilon_{n}$ in formula
\eqref{2.6}. We obtain
\[ (C^{2}+ C_{2})(\|\theta(x,E+i\epsilon_{n})\|_{L_{2}(B_{R_{n}}\cap \Omega)}+ \|\phi\|_{L^2})
\|u\|_{L_{2}(B_{R_{n}}\cap \Omega)} \geq |c|(R_{n}-R_{0}-1). \]
Substituting \eqref{2.7}
into the last inequality, we find that there exists some
constant $C_{3}$ such that for $n$ large enough, we have
\begin{equation}
\label{2.8}
\|\theta(x,E+i\epsilon_{n})\|_{L_{2}(B_{R_{n}}\cap \Omega)}
\geq C_{3}R_{n}^{1-\frac{\alpha}{2}}-\|\phi\|_{L^2}. \end{equation}
Now it remains to invoke Lemma 2.1 and note that
\[ \Im ((H_{V}-E-i\epsilon_{n})^{-1}\phi, \phi) \geq \epsilon_{n}
\|\theta(x,E+i\epsilon_{n})\|^{2}_{L_{2}(B_{R_{n}})} \]
for every $n.$ Using the estimate \eqref{2.8} and the relation
between $R_{n}$ and $\epsilon_{n},$ we find
\[ \Im ((H_{V}-E-i\epsilon_{n})^{-1}\phi, \phi) \geq C_{4}
\epsilon_{n}^{\alpha -1} \]
for sufficiently small $\epsilon_n$
The application of Lemma 2.2 now completes the proof.
\end{proof}  \rm

\noindent \it Remarks. \rm 1. Theorem 1.1 also holds for wider classes
of potentials and boundary conditions. 
The restrictions of the classes
come from Lemma 2.3, the necessary estimate on the energy norms.
With the help of smooth mollifiers to justify integration by parts,
Theorem 1.1 can be extended to the classes to which
one can extend Lemma 2.3.

\smallskip
2. We also note that the same argument as in the proof implies
that $D^{\alpha}\mu^{\phi}(E) = \infty$ if instead of \eqref{1.2}
in the assumption of Theorem 1.2 we suppose that
\[ \liminf\limits_{R\rightarrow \infty} R^{-\alpha}\|u(x,E)\|^{2}_{B_{R}}=0. \]
We will use this fact in the proof of Corollary 2.6 below. \\

The next question that we would like to discuss is a sufficient condition
for the existence of the various components of the spectrum.
Let us recall the definition of Hausdorff measures and dimension. For
$\alpha \in [0,1]$ and any $S \subset \R,$ the $\alpha$-dimensional
Hausdorff measure of $S$ is defined by
\[ h^{\alpha}(S) = \lim\limits_{\delta \rightarrow 0} \inf_{\delta -
\rm{covers}} \sum\limits_{\gamma =1}^{\infty} |I_{\gamma}|^{\alpha}, \]
where $I_{\gamma}$ are the intervals constituting the cover.
The Hausdorff dimension of a set $S$ is the infimum of all values of $\alpha$
such that $h^{\alpha}(S)=0.$
 First, we are going to prove \\

\gap
\bf Theorem 2.5. \it Let $H^{\Omega}_{V}$ be a Schr\"odinger
 operator, with $V$ and $\Omega$ satisfying the same conditions
as in Theorem 1.1.
Suppose that for a measurable set $S$ of
positive $h^{\alpha}$ measure, for each $E \in S$, there
exists a non-trivial solution $u(x, E)$ of the generalized
eigenfunction equation \eqref{1.1} satisfying the boundary
conditions such that
\[ \liminf_{R \rightarrow \infty} R^{-\alpha}
\|u(x, E)\|^{2}_{B_{R}} < \infty. \]
Then there exists a vector $\varphi \in L^{2}(R^{n})$ such that
$\mu^{\varphi}(S_{1})>0$ for any $S_{1} \subset S$ of positive
$h^{\alpha}$ measure. In particular, if $\alpha=1,$
we have absolutely continuous spectrum filling the set $S.$

\gap
\noindent Remark. \rm In many applications, particularly in one
dimension, one applies a reasoning different from the one suggested by Theorem
2.5 to derive existence of various dimensional spectral components
from results like Theorem 1.1.
One proves the existence of solutions as in (2) for a.e.~$E$, and then
uses rank-one perturbation arguments (see, e.g., \cite{JiLai, KLS}). 

\begin{proof} \rm 
Recall that for every self-adjoint operator there is an associated 
spectral measure of maximal type, $\mu,$ such that for every $\psi$ and 
any measurable set $S,$ $\mu^{\psi}(S)>0$ implies $\mu(S)>0.$ A vector 
$\chi$ is of the maximal type if for any measurable
set $S,$ $\mu^{\chi}(S)>0$ given that $\mu(S)>0.$ We will show that for any
$S_{1} \subset S$ of positive $\alpha$-dimensional Hausdorff measure,
there exists
a vector $\psi$ with $\mu^{\psi}(S_{1})>0.$
By the standard argument for the existence of vectors of maximal
type (see e.g. \cite{BS}), this would imply existence of the vector
$\varphi$ as in the theorem.
Pick some ball $B_{R_{0}}$
such that $\|u(x,E)\|_{L^{2}(B_{R_{0}\cap \Omega})} \neq 0$
for energies $E$ in a subset $S_{2}$ of $S_{1}$ of positive
$h^{\alpha}$ measure (it is easy to see that such a ball exists,
because of the $\sigma$-additivity
of $h^{\alpha}$). We remark that for a wide class of operators
$H_{V}^{\Omega},$ an arbitrary ball will do because of the unique
continuation (solutions $u(x, E)$ cannot vanish identically on
any ball), but there is no need to invoke these results.
Pick a basis $\{ \psi_{n}(x) \}_{n=1}^{\infty}$ in the Hilbert space
$L^{2}(B_{R_{0}} \cap \Omega)$. Since $\{\psi_{n}\}$ forms a basis, for
every $E \in S_{2}$ there exists an $n$ such that
\[ \int\limits_{B_{R_{0}}\cap \Omega} \psi_{n}(x)u(x,E)\,dx
\neq 0. \]
Consider the functions $D^{\alpha}\mu^{\psi_{n}}$
on the set $S_{2}.$ By Theorem 1.1, for every $E \in S_{2}$
there exists an $n$ such that
$D^{\alpha}\mu^{\psi_{n}}(E)>0.$
In particular, by $\sigma$-additivity of $h^{\alpha}$, there exists
an $n_{0}$ such that
$D^{\alpha}\mu^{\psi_{n_{0}}}(E)>0$  for
every $E$ in a set $S_{n_{0}} \subset S_{2}$ of positive $h^{\alpha}$
measure.
By the results of Rogers-Taylor theory (see \cite{Ro}, Theorem 63), it follows
that the measure $\mu^{\psi_{n_{0}}}$ gives positive
weight to the set $S_{n_{0}},$ and hence to the set $S_{1}.$  The case
of the absolutely continuous spectrum corresponds to $\alpha=1;$
in this case the application of Rogers-Taylor theory may be replaced
by a well-known fact that
a measure gives positive weight to a set of positive
Lebesgue measure when its derivative is positive a.e.\ in this set.
\end{proof} \rm

\rm From Theorem 2.5 (or, essentially, from its proof and the remark
after the proof of Theorem 1.1) immediately
follows:  \\

\gap
\bf Corollary 2.6. \it For any $\alpha,$ the set $S$ of
energies
$E$ for which there exists a solution $u(x,E)$ satisfying
 \begin{equation}
\label{2.9}
\liminf_{R \rightarrow \infty} R^{-\alpha}
\|u(x, E)\|^{2}_{B_{R}} =0 \end{equation}
has zero $h^{\alpha}$ measure.  \\
\gap

\noindent \it Remark. \rm The fact that there may be only
countably many values
of $E$ (counting multiplicities) for which equation
\eqref{1.1} has
$L^2$ solutions satisfying the boundary conditions,
is an obvious consequence of
the separability of the Hilbert space $L^{2}(\Omega).$ This corollary
may be viewed as a  less trivial generalization for slower rates
of decay. 

\begin{proof} \rm Suppose that $S$ has positive $h^{\alpha}$ measure.
By the remark after the proof of Theorem 1.1, \eqref{2.9} implies that
$D^{\alpha}\mu^{\phi}(E)=\infty$ for every $E \in S$ and finitely
supported $\phi$ such that $\int u(x,E)\phi(x) \ne 0.$
 Proceeding as in
the proof of Theorem 2.5, we can find a vector $\varphi$ such that
$D^{\alpha}\mu^{\varphi}(E)=\infty$ for any $E$ in some
 set of positive $h_{\alpha}$ measure. This is not possible by Rogers-Taylor 
(see \cite{Ro}, Theorem 67) and therefore gives a contradiction.
We remark that for $\alpha =1$, this argument reduces to the well-known
statement
that a finite Borel measure $\mu^{\psi}$ cannot have an infinite derivative
on a set of positive Lebesgue measure.
\end{proof} \rm

We would like to end this section by drawing a link with the well-known
results of Rellich \cite{Rel} and Kato \cite{Ka} who showed, respectively,
that for the free Laplacian and the Laplacian with a short-range perturbation
(i.e., a potential which satisfies $|V(x)| \leq C(1+|x|)^{-1-\epsilon}$),
there are no solutions satisfying \eqref{2.9}
 with $\alpha =1$ for \it any energy. \rm
Corollary 2.6 shows that for a much larger class of potentials, such
solutions are still in some sense ``exceptional'' and can only occur
on a set of energies of zero Lebesgue measure.

\gap
\begin{center} \bf \large 3. Solutions and spectrum: Discrete case
\end{center}
\mgap
In this section, we consider discrete Schr\"odinger operators.
All the results of the previous section extend to the discrete setting.
In fact, the proofs are easier due to the absence of the Sobolev
estimates issue, and there are no restrictions on potential.

Let $\Omega$ be some connected infinite domain in $\Z^{d}.$
We define the Schr\"odinger operator $h_{v}^{\Omega}$ on
$L^{2}(\Omega)$ with Dirichlet
boundary conditions by
\[ h_{v}^{\Omega}f(n) = \sum\limits_{|n-m|=1,\, m \in \Omega}f(m)+
v(n)f(n). \]
It is easy  to check that the operator defined in this way is self-adjoint.
 
We need an analog of the Green's formula in the discrete setting.
For any domain $S \subset Z^d$ let us denote by $\partial S$ the set of 
points outside $S$ which have a point of $S$ within a unit distance.
 We have for any two functions
$f,$ $g$ 
\[
 \sum\limits_{n \in S}(h^{\Omega}_{v}f(n)g(n) -
f(n) h^{\Omega}_{v}g(n)) 
 =\sum\limits_{m \in \partial S} \left(
f(m)\sum\limits_{l \in N_{S}(m)} g(l)
 - g(m)\sum\limits_{l \in N_{S}(m)} f(l)\right), 
\]
where $N_{S}(m)$ denotes the set of neighbors of the point $m \in 
\partial S$ lying 
in $S$ (so that $|m-n|=1$ for any $n \in N_{S}(m)$).
 We will say therefore that the analog of the 
Wronskian over $\partial S$ of two functions is,
in the discrete setting,
\[ w_{\partial S}[f,g] = \sum\limits_{m \in \partial S} \left(
f(m)\sum\limits_{l \in N_{S}(m)} g(l)
 - g(m)\sum\limits_{l \in N_{S}(m)} f(l) \right). \]

For convenience, in all considerations for the discrete case,
we replace the balls $B_{R}$ with cubes $C_{R}.$
 The point $n=(n_{1}, \dots n_{d})$ of the lattice belongs
to $C_{R}$ if and only if  $|n_{i}| \leq R$ for all $i=1, \dots, d.$
 
 We now formulate and prove an analog of Theorem 1.1 in the discrete case. \\

\gap
\bf Theorem 3.1. \it
Suppose that there exists a solution $u(n, E)$ of the generalized
eigenfunction equation
\begin{equation}
\label{3.1}
 (h_{v}^{\Omega}-E)u(n,E)=0
\end{equation}
satisfying the Dirichlet  boundary conditions on $\partial \Omega.$ 
 Suppose that for
some $\alpha,$ $0 \leq \alpha \leq 1,$ we have
\begin{equation}
\label{3.2}
 \liminf_{R \rightarrow \infty}
R^{-\alpha}\sum\limits_{n \in C_{R} \cap \Omega}|
u(n,E)|^{2}\,dx < \infty.
\end{equation}
Fix some vector $\phi$ of compact support such that
\[ \sum\limits_{n} u(n, E)\phi(n) \ne 0. \]
Then we have
\[ D^{\alpha}\mu^{\phi}(E) >0. \]
In particular, if $u(n_{0},E) \ne 0,$ then
\[   D^{\alpha}\mu^{\delta_{n_{0}}}(E) >0. \]
{\rm(}here $\delta_{n_{0}}$ is a function equal to $1$ at $n_{0}$ and
$0$ otherwise{\rm)}. \\
\gap
\begin{proof}
 The argument repeats the proof of Theorem 1.1, except that we
do not need Lemma 2.3. 
The analog of Lemma 2.4 is proven directly by 
the observation
that 
\[ w_{\partial(\Omega \cap C_{r})}[f,g]=
w_{\partial C_{r} \setminus \partial \Omega}[f,g], \]
and 
\[ \sum\limits_{r=1}^{R} |w_{\partial C_{r} \setminus \partial \Omega}
[f,g]|
  \leq d 
\|f\|_{L^{2}(C_{R+1}\cap \Omega)} \|g\|_{L^{2}(C_{R+1}\cap \Omega)}. \]
\end{proof} \rm

\gap
\it Remark. \rm As in the continuous case, the result is also true for 
more general boundary conditions. \\

We also have an analog of Theorem 2.5: \\

\gap
\bf Theorem 3.2. \it Suppose that for each energy $E$ in some
 measurable set $S$ of positive
$h^{\alpha}$  measure, there exists a non-trivial solution $u(n, E)$
of equation \eqref{3.1}
satisfying Dirichlet boundary conditions and having the property \eqref{3.2}.
Then there exists a vector $\phi \in L^{2}(\Z^{d})$, such that
$\mu^{\phi}(S_{1})>0$ for any set $S_{1} \subset S$
of positive $h^{\alpha}$ measure. In particular, if $\alpha=1,$
the set $S$ is an essential support of the absolutely continuous part
of the measure
$\mu^{\phi}$ restricted to $S.$ \\

\gap
\noindent \rm The proof of this theorem is the same as the proof of
Theorem 2.5.

\gap
\begin{center} \bf \large 4.  Examples and discussion \end{center}
\mgap
The purpose of this
section is purely illustrative --- to show where the solutions
we are studying are known to occur. However, these observations
will also partly
lead us to the issue which is the topic of the next section:
the relationships between generalized eigenfunctions, spectrum, and
dynamics.
In addition, we show that  the criteria
given by Theorems 1.1 and 3.1 are sufficient but not, in general, necessary
for the positivity of the derivatives of spectral measures. We give
an explicit example to confirm this statement.

Our first remark is that solutions $u(x,E)$ satisfying
\[ \liminf_{R \rightarrow \infty}R^{-1} \|u(x,E)\|^{2}_{B_{R}} < \infty \]
exist for every energy $E \ne 0$
in the spectrum in the case of the free Laplacian operator
in $\R^{d}$
or in the cylinder with Dirichlet boundary conditions.
In the cylinder case, we may take
\[ u(x,E)= \exp( i \sqrt{E-E_{l}}x_{1}) Z_{\overline{l}}
(x_{2}, \dots x_{d}), \]
where $x_{1}$ is the coordinate along the rotation axis, $E_{l}$ is any
eigenvalue (less then $E$) of the Laplace operator with Dirichlet boundary
conditions on the
$d-1$-dimensional
ball, and $Z_{\overline{l}}
(x_{2}, \dots, x_{d})$ is any eigenfunction corresponding
to this eigenvalue.
 In the free
case, we can take any function
\[ u(x,E) = r^{-\frac{d}{2}+1} J_{\nu}(\sqrt{E}r)Y_{\overline{l}}
(\overline{\theta}), \]
 where  $Y_{\overline{l}}$ is
any of the spherical harmonics corresponding to the eigenvalue
$E_{l}=l(l+d-2)$ of the Laplace-Beltrami operator on the $d$-dimensional
sphere,
and $J_{\nu}$ is a Bessel function (without singularity at the
origin) with $\nu$ defined by $\nu^{2}= l(l+d-2)+(\frac{d}{2}-1)^{2}.$
Note that for large $r,$
\[ J_{\nu}(\sqrt{E}r) \sim Cr^{-\frac{1}{2}}\cos\left(\sqrt{E}r-
\frac{2\pi \nu- \pi}{4}\right)(1+o(1)). \]
See, for example, \cite{Erd,Vil} for more information on
spherical harmonics and Bessel functions.

Using the results of Agmon theory and related estimates on the
Fourier transform (see \cite{Ag} or \cite{ReSi}, and \cite{AH}), it is 
straightforward to show
that the existence for every $E \in (0,\infty)$
of solutions with the rate of growth of the $L^{2}$
norm as in \eqref{1.2} with $\alpha=1$ extends to
perturbations of the free Laplacian by
short range potentials, $|V(x)| \leq
C(1+|x|)^{-1-\epsilon},$ if $C$ is sufficiently small.
In one dimension, it was recently shown \cite{CK,Re} that such solutions
exist for a.e.\ $E \in (0,\infty)$ for any potential $V$ satisfying
$|V(x)| \leq C(1+|x|)^{-\frac{1}{2}-\epsilon}.$ This implies that the
absolutely continuous spectrum of the free operator in one dimension
is stable under all perturbations decaying at this rate.
This result is optimal: there are potentials which satisfy
$|V(x)| \leq C(1+|x|)^{-\frac{1}{2}}$ and lead to purely singular
spectrum in $(0, \infty).$ The corresponding question about
the borderline decay for the stability of the absolutely continuous
spectrum is open in higher dimensions, with any power in $[1,\frac{1}{2}]$
a possible candidate, in principle.
We conjecture \\

\gap
\bf Conjecture I. \it Suppose that $H_{V}$ is a Schr\"odinger operator in $\R^{d}$
for which $|V(x)| \leq C(1+|x|)^{-\frac{1}{2}-\epsilon},$ $\epsilon >0.$
Then the absolutely continuous spectrum of the operator $H_{V}$
fills the whole positive
semi-axis. \\
\gap

\rm This conjecture would in particular follow from \\

\gap
\bf Conjecture II. \it Under the conditions of the previous conjecture, for
a.e.\ $E \in (0,\infty)$ there exists a solution $u(x,E)$ of the
generalized eigenfunction equation satisfying
\eqref{1.2} with $\alpha=1.$ \\
\gap

\rm Our next  example concerns Schr\"odinger operators with
periodic potentials. Let $V(x)$ be a smooth periodic potential of period
one in all variables $x_{1}, \dots, x_{d}.$
Given $E$ in the spectrum of $H_{V},$ consider the boundary value
problem
\begin{eqnarray}
 && (H_{V}-E)b(x,E) =0,  \nonumber \\
&& \frac{\partial^{j}b}{\partial x_{l}^{j}} \left|_{x_{l}=1} \right.
= \exp(i \theta_{l}) \frac{\partial^{j}b}{\partial x_{l}^{j}}
\left|_{x_{l}=0}, \right. \,\,\, l=1, \dots d, \,\,\, j=0,1. 
\label{4.1}
\end{eqnarray}
The set of all values of $\theta \in [0,2\pi)^{d}$
for which there exist solutions of the boundary value problem \eqref{4.1}
is called the real (physical) Fermi surface $F_{E}.$ From well-known
results on spectral properties of periodic differential operators
(see \cite{Kuch}), it follows that for all but a countable
set of energies in the spectrum (exceptional points corresponding
to band edges), we can find solutions $u(x,E)$ of the generalized
eigenfunction equation \eqref{1.1} of the following type:
\begin{equation}
\label{4.2}
 u(x,E) = \int\limits_{S} b(x,\theta, E)
\gamma (\theta ) d\sigma,
\end{equation}
where $S \subset F_{E}$ is a piece of
an analytic $(d-1)$-dimensional surface, $\gamma(\theta)$
is a $C_{0}^{\infty}(S)$-function and $b(x,\theta,E)$
are Bloch functions satisfying \eqref{4.1}
\[ b(x, \theta, E) = \exp (i \theta x)
  f(x, \theta,E), \]
where $ f(x,\theta, E)$ is periodic with period one
in all directions in $x,$ continuous in $x,$ and
analytic (as an $L^{2}([0,1)^{d})$ vector) in $\theta \in S.$
We claim that $u(x,E)$ satisfies
\[ \liminf_{R \rightarrow \infty} R^{-1}\|u(x,E)\|^{2}_{B_{R}} \leq \infty. \]
This can be shown in a way similar to the proof of this property in the case of
Fourier transforms of measures supported on $(d-1)$-dimensional
smooth surfaces (see \cite{AH}).
Represent the equation of the surface $S$ as $\theta_{d} = s(\theta_{1}
, \dots, \theta_{d-1})$ (we can assume that $S$ is small enough and
$\theta_{d}$ is chosen so that this is possible). Then we can
rewrite \eqref{4.2} as
\[ u(x,E) = \int\limits_{S'} \exp (i \theta'x' + i s(\theta')x_{d})
f(x, \theta',E) \gamma'(\theta') d\theta', \]
where the integration is now over the projection $S'$ of $S$ on
the hyperplane $\theta_{d}=0,$
 $\theta'$ denotes first $d-1$ coordinates,
and $\gamma'$ includes the Jacobian from the change of variables.
 Fix the value of $x_{d}$ and integrate over the cube $C_{R}'$ in the other
coordinates $x'=x_{1}, \dots, x_{d-1}:$  
\begin{eqnarray*}
 \int\limits_{C_{R}'} |u(x,E)|^{2} dx' & = &
\int\limits_{S'}\int\limits_{S'} \gamma'(\theta') \gamma'(\tilde{\theta}')
\exp (i(s(\theta')-s(\tilde{\theta}'))x_d)  \\
& &  \int\limits_{C'_{R}}
\exp(i (\theta'-\tilde{\theta}')x') f(x', \theta',E) f(x', \tilde{\theta}',E)
\,dx' d \theta' d\tilde{\theta}'. 
\end{eqnarray*}
Without loss of generality, take $R$ to be an integer. Then
we obtain
\begin{equation}
\label{4.3}
\int\limits_{C_{R}'} |u(x,E)|^{2} dx'=  \int\limits_{S'}\int\limits_{S'}
 d \theta' d\tilde{\theta}'
\left(\prod\limits_{j=1}^{d-1}\frac{\sin (R+\frac{1}{2})
(\theta_{j}-\tilde{\theta}'_{j})}{\sin \frac{1}{2}( \theta_{j}
- \tilde{\theta}'_{j})} \right) \psi (\theta', \tilde{\theta}'),
\end{equation}
where
\begin{eqnarray*}
 \psi (\theta', \tilde{\theta}')& = & \gamma'(\theta') \gamma'(\tilde{\theta}')
\exp (i(s(\theta')-s(\tilde{\theta}'))x_d) \\
& & \int\limits_{C_{1}'} f(x', \theta',E) 
f(x', \tilde{\theta}',E)
\exp(i (\theta'-\tilde{\theta}')x')\,dx'. 
\end{eqnarray*}
Due to the properties of $f$ and $\gamma,$ the function $\psi$ is smooth
and hence the right-hand side in \eqref{4.3}
 converges as $R \rightarrow \infty$
to the constant
\[ C= \int\limits_{S'}
 d \theta'  \psi(\theta', \theta') = \int\limits_{S'} d\theta'
|\gamma'(\theta')|^{2} \left( \int\limits_{C_{1}'} |f(x', \theta'),E|^{2}\,dx'
\right). \]
Therefore integrating in $x_{d}$ from $-R$ to $R,$ we obtain
\[  \int\limits_{B_{R}} |u(x,E)|^{2} dx \leq C R, \]
as claimed.

\begin{figure}
\begin{center}
\includegraphics{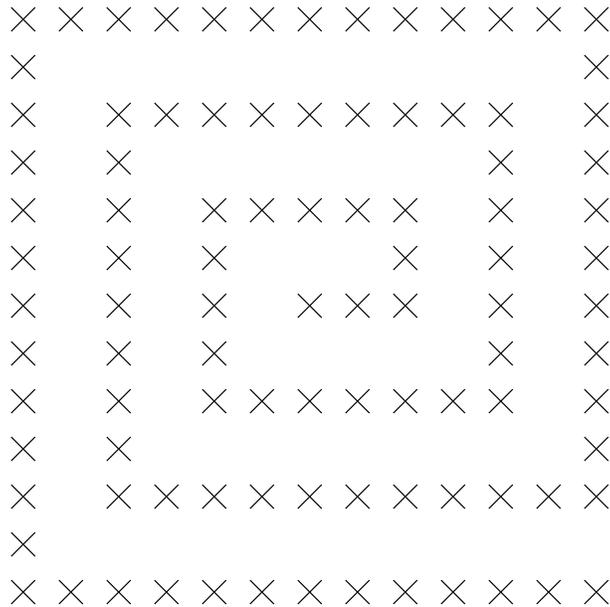}
\caption{The spiral domain}\label{Fi:A}
\end{center}
\end{figure}

Our last example in this section shows that the criteria for
the positivity of derivatives of spectral measures,
given by Theorems 1.1 and 3.1, provide
 a sufficient, but, in general, not necessary condition. The example is
especially simple and transparent
in the discrete setting. Let us consider the discrete plane $\Z^{2}$
and let $\Omega$ be an infinite ``spiral'' in this plane (see Figure 1;
we marked by $\times$ the points which do not belong  to the domain).
Consider $h_{0}^{\Omega}$ defined on the spiral with Dirichlet boundary
conditions. By inspection, we see that $h_{0}^{\Omega}$ acts on
$l^{2}(\Omega)$ as a free one-dimensional Jacobi matrix. Hence the
spectrum is absolutely continuous in $[-2,2],$ and for every $E$ in
this interval, there exists an explicitly computable unique solution
$u(n,E)$ of
the generalized eigenfunction equation satisfying the boundary
conditions:
\[ u(n,E) = \sin \left(\cos^{-1}\left(\frac{E}{2}\right) n\right). \]
 This is a standard discrete plane wave. If we measure the
linear distance $N$ along the spiral, the square of the $l^{2}$ norm of
this solution grows as $N.$ However in $\Z^{2},$ we have
\[ \|u(x,E) \|^{2}_{B_{R} \cap \Omega} \sim R^{2}. \]
Hence in this case, we cannot find solutions as in Theorem 1.1.

 We remark that \cite{Sim2} contains an example
of a \it bounded \rm spiral ``jelly roll'' domain on which the Laplace
operator with Neumann boundary conditions has
 absolutely continuous spectrum. In this case, for a.e.\
$E$ in the spectrum, the norm of solutions becomes infinite for
finite $R.$

\gap
\begin{center} \bf \large 5. Solutions and dynamics \end{center}
\mgap
In this section, we prove Theorem 1.2 and apply it to study quantum
dynamics in the random decaying potentials model studied
in \cite{DSS} and more recently in \cite{KLS,KU}.

The previous section provided us with several examples of operators
with absolutely continuous spectrum and solutions satisfying the
condition \eqref{1.2}
in Theorem 1.1 for $\alpha=1,$ and one example of an operator with
absolutely continuous spectrum, but without such solutions.
For the former three, the transport is ballistic for every
vector (i.e., $\langle \langle |X|^{m} \rangle \rangle_{T}
\sim T^{m}$); for the latter, it is easy to see that the transport
is not ballistic (it is diffusive in $\Z^{2}$). Theorem 1.2 indicates
that this is not a coincidence.

The proof of Theorem 1.2 is an extension of the proof of Theorem 6.1
of \cite{L}, and it is essentially the same in both the
discrete and continuous settings. We will use a discrete notation
which formally only covers the discrete case, but the continuous
case follows from it in a totally straightforward manner (which
essentially amounts to replacing $n$ by $x$ and some summations
by integrals). We note that in \cite{L} the continuous case
(Theorem 6.2 of \cite{L}) is getting an independent treatment, based
on semi-group kernel inequalities. This is not needed here, since
we assume the existence of eigenfunction expansions with suitable
properties. This allows our Theorem 1.2 to cover some cases, such
as Stark operators, that are excluded from Theorem 6.2 of \cite{L}.

Recall that a measure $\mu$ is called uniformly $\alpha$-H\"older
continuous (denoted U$\alpha$H)
if there exists a constant $C$ such that for every interval $I$
with $|I|  \leq 1$ we have
\begin{equation}
\label{5.1}
\mu(I) \leq  C|I|^{\alpha}.
\end{equation}
$\alpha$-continuous measures (recall that this means
 measures giving zero weight to all sets of zero $h^{\alpha}$
measure) can be approximated by U$\alpha$H measures in the
following sense: \\

\gap
\bf Theorem (Rogers-Taylor {\rm \cite{RT}}). \it
A finite Borel measure $\mu$ on $\R$ is $\alpha$-continuous if
and only if for every
$\epsilon >0,$ there exist two mutually singular Borel measures
$\mu_{1}^{\epsilon}$ and $\mu_{2}^{\epsilon}$, such that
$\mu = \mu_{1}^{\epsilon} + \mu_{2}^{\epsilon}$, where
$\mu_{1}^{\epsilon}$ is U$\alpha$H and $\mu_{2}^{\epsilon}(\R)<
\epsilon.$ \\

\gap
\rm For U$\alpha$H measures, we can study dynamics
with the aid of the following Strichartz estimate: \\

\gap
\bf Theorem (Strichartz {\rm \cite{Str}}). \it
Let $\mu$ be a finite
U$\alpha$H measure, and for each $f \in L^{2}(\R, d\mu)$ denote
$$ \widehat{f\mu}(t)  = \int \exp(-ixt)f(x) \,d\mu(x). $$
Then there exists a constant $C_{1},$ depending only on $\mu$
{\rm(}more precisely, only on $C$ in \eqref{5.1}{\rm)}, such that for any
$f \in L^{2}(\R, d\mu)$ and $T>0$
$$ \langle |\widehat{f\mu}|^2 \rangle_{T} < C_{1} \|f\|^{2} T^{-\alpha}, $$
where $\|f\|$ is the $L^{2}$ norm of $f$.  \\

\gap
\rm We now prove Theorem 1.2.\\
\begin{proof}
\rm
Without loss of generality, assume $\|\psi\|=1.$
We first establish the existence of a Borel set $\tilde S \subset S$,
for which the following three properties are true: \\
(i) $\|P_{\tilde S}\psi\| > 0$. \\
(ii) The restriction of the spectral measure $\mu^\psi$ to $\tilde S$
is U$\alpha$H. \\
(iii) There exists a constant $C_2$, such that for each $E\in \tilde S$
and $R > 0$,
the corresponding generalized eigenfunction $u(n,E)$ satisfies
$$\sum_{|n|< R} |u(n,E)|^2 < C_2 R^\gamma .$$

\smallskip
We shall establish (i)--(iii) in two stages. First, note that the function
$$f(E)\equiv\sup_{R > 0} R^{-\gamma}\sum_{|n|< R} |u(n,E)|^2$$ is a measurable
function of $E$, which, by \eqref{1.3}, is finite everywhere on
$S$. Thus, since $S = \bigcup_{k=1}^{\infty} \{E\in S\,|\, f(E) < k\}$,
there is clearly a Borel subset $S_1\subset S$ of positive $\mu^\psi$
measure and a constant $C_2$, such that $f(E) < C_2$ for any $E\in S_1$.
That is, property (iii) holds for $S_1$. Next, since the restriction
of $\mu^\psi$ to $S_1$ is $\alpha$-continuous, the above mentioned
Rogers-Taylor theorem implies that there is a Borel subset 
$\tilde S\subset S_1$ of positive $\mu^\psi$ measure (so property (i)
holds) such that the restriction of $\mu^\psi$ to $\tilde S$ is
U$\alpha$H (so property (ii) holds).

Let us now denote $\psi_1 = P_{\tilde S}\psi$,
$\psi_2 = P_{\R\setminus\tilde S}\psi$, where $P_{\cdot}$ denotes the
spectral projection over the corresponding set. Then $\psi = \psi_1 + \psi_2$,
and $\psi_1$, $\psi_2$ are mutually orthogonal so
$1 = \|\psi\|^2 = \|\psi_1\|^2 + \|\psi_2\|^2$.
Let $P_{R_{T}}$ be the projector on the set of sites $n$ with $|n| \leq R_{T}$.
$R_{T}$ is a function of the time parameter $T$ to be chosen later. 
Given any vector $\varphi$, we will routinely use the notation
$\varphi(t)= \exp(-ih_{v}^{\Omega}t) \varphi$.

We have 
\begin{eqnarray*}
\langle \|P_{R_{T}}\psi_{1}(t)\|^{2}\rangle_{T}  & = &
\sum\limits_{|n| \leq R_{T}} \frac{1}{T} \int\limits_{0}^{T}
\left|\int \exp(-iEt) \overline{u(n,E)}\, d \mu^{\psi_{1}}(E)\right|^2 \,dt \\
& \leq & C_{1} T^{-\alpha}\sum_{|n|< R_{T}} \int |u(n,E)|^{2} d\mu^{\psi_{1}}(E) \\
& \leq & C_{1} 
\|\psi_{1}\|^{2} \left( \sup_{E\in{\tilde S}} \sum_{|n|< R_{T}} |u(n,E)|^{2} \right)
T^{-\alpha}
\end{eqnarray*}
by Strichartz theorem, and so
\begin{equation}
\label{5.2}
\langle \|P_{R_{T}}\psi_{1}(t)\|^{2}\rangle_{T} \leq
C_2 C_1 \|\psi_1\|^2 R_T^\gamma T^{-\alpha} .
\end{equation}

For each $T>0$, we now define
$$R_T = \left({\frac{\|\psi_1\|^2\,T^\alpha} 
{64\,C_2 C_1}}\right)^{1/\gamma} ,$$
such that we have 
$$\langle\|P_{R_T}\psi_1(t)\|^2\rangle_T <
{\frac{\|\psi_1\|^4}{64}} ,$$
and thus
\begin{eqnarray*}
\langle\|P_{R_T}\psi(t)\|^2\rangle_T &\leq &
\left\langle\left(\|P_{R_T}\psi_1(t)\|+\|P_{R_T}\psi_2(t)\|
\right)^2\right\rangle_T \\
&\leq& \left\langle\left(\|P_{R_T}\psi_1(t)\|+\|\psi_2\|
\right)^2\right\rangle_T \\   
&\leq& \left(\sqrt{\langle\|P_{R_T}\psi_1(t)\|^2\rangle_T}+\|\psi_2\|
\right)^2 \\   
&<& \left({\frac{\|\psi_1\|^2}{8}}+\|\psi_2\|\right)^2 \\
&=& {\frac{\|\psi_1\|^4}{64}}+\|\psi_2\|^2+{\frac{1}{4}} 
\|\psi_2\| \|\psi_1\|^2 \\  
&<& \|\psi_2\|^2+{\frac{1}{2}}\|\psi_1\|^2 \\
&=& 1-{\frac{1}{2}}\|\psi_1\|^2 .
\end{eqnarray*}
Since
$$\langle\|P_{R_T}\psi(t)\|^2\rangle_T + 
\langle\|(1-P_{R_T})\psi(t)\|^2\rangle_T = 1 ,$$
we obtain
$$\langle\|(1-P_{R_T})\psi(t)\|^2\rangle_T > 
{\frac{1}{2}}\|\psi_1\|^2 ,$$
which implies
$$\langle\langle|X|^m\rangle\rangle_T >  
{\frac{1}{2}}\|\psi_1\|^2 R_T^m = 
{\frac{\|\psi_1\|^2}{2}} \left({\frac{\|\psi_1\|^2} 
{64\,C_{2}C_1}}\right)^{m/\gamma}\,T^{\alpha m/\gamma} ,$$
proving \eqref{1.4}.
\end{proof} \rm

Note that the above proof does not attempt to provide optimal
estimates. We could (by allowing various constants to grow) choose
$\psi_1$ to have a norm that is arbitrarily close to that of
$P_S\psi$, and $R_T\sim T^{\alpha/\gamma}$ so that 
$\langle\|(1-P_{R_T})\psi(t)\|^2\rangle_T$ is larger than
something arbitrarily close to $\|P_S\psi\|^2$. This means
that there is a component of the wave packet of size corresponding
to $\|P_S\psi\|$ that is spreading on average at a rate of at
least $T^{\alpha/\gamma}$.

We now apply Theorem 1.2 to investigate dynamics for the following
model. Let $v_{\omega}(n)$ be independent  random variables
such that
\begin{equation}
\label{5.3}
E(v_{\omega}(n)) =0, \,\,\,\,\,E(v_{\omega}(n)^{2})^{\frac{1}{2}}
= \lambda n^{-\frac{1}{2}},\,\,\, {\rm and}
\,\,\,\,\sup_{\omega}|v_{\omega}(n)|
\leq Cn^{-\frac{1}{3}-\delta}, \,\,\,\delta>0.
\end{equation}
 For example, if we take
i.i.d.  random variables $a_{\omega}(n)$
with uniform distribution in $[-\sqrt{3},\sqrt{3}]$,
then $v_{\omega}(n) = \lambda n^{-\frac{1}{2}} a_{\omega}(n)$
satisfy all the conditions. The half-line random Schr\"odinger operators
$h_{\omega}$ with, say, Dirichlet boundary conditions at zero and
potential $v_{\omega}$  exhibit very rich
spectral structure. Such operators where studied by Delyon, Simon, and
Souillard \cite{DSS}, and more recently by Kotani
and Ushiroya  \cite{KU}, and by Kiselev, Last, and Simon \cite{KLS}.
Our study here is based mainly on the results of the last
paper. In particular, the following has been proven in \cite{KLS}: \\

\gap
\bf Theorem (KLS {\rm \cite{KLS}}). \it
For all $\omega$, the essential spectrum of $h_{\omega}$ is $[-2,2].$
If $|\lambda| < 2,$ then for a.e.\ $\omega$, $h_{\omega}$ has purely
singular continuous spectrum in
$\{ E \,| \,|E| < (4 -\lambda^{2})^{\frac{1}{2}} \}$
and only dense pure point spectrum in
$\{ E \,|\, (4 -\lambda^{2})^{\frac{1}{2}}< |E| < 2  \}.$

For a.e.\ $\omega$ and $E \in (-2,2)$
\begin{equation}
\label{5.4}
\lim\limits_{n \rightarrow \infty} \frac{\log \|T_{E}(n,0)\|}
{\log n}=  \frac{\lambda^{2}}{8-2E^{2}},
\end{equation}
and there
exists an initial condition $\theta(\omega)$ at zero such that
\begin{equation}
\label{5.5}
\lim\limits_{n \rightarrow \infty} \frac{\log \|T_{E}(n,0)u_{\theta(\omega)}\|}
{\log n}=  -\frac{\lambda^{2}}{8-2E^{2}},
\end{equation}
where $u_{\theta(\omega)}$ is the $2$-vector corresponding to the
boundary condition $\theta(\omega)$ at $0,$ and
$T_{E}(n,0)$ is the transfer matrix from $0$ to $n$ at energy $E$. \\

\gap
\rm This theorem implies that for a.e.\ $\omega$ and
$E\in (-\sqrt{4 - \lambda^{2}}, \sqrt{4 -\lambda^{2}})$, the spectral
measure $\mu$ (corresponding to the vector $\delta_{1}$) has
local Hausdorff dimension
\begin{equation}
\label{5.6}
 \alpha(E, \lambda) = \frac{4-E^{2}-\lambda^{2}}{4-E^{2}}
\end{equation}
at energy $E$, in the sense that for any $\epsilon > 0,$ there is
a $\delta$ so that $\mu(A)=0$ if $A$ is a subset of $(E-\delta,
E+\delta)$ of Hausdorff dimension less than $\alpha(E,\lambda)-\epsilon$,
and there is a subset $B$ of Hausdorff dimension less than
$\alpha(E, \lambda)+\epsilon$ such that $\mu((E-\delta, E+\delta) \setminus B) = 0$.
These properties of the spectral measure follow from \eqref{5.4}, \eqref{5.5}
by subordinacy theory \cite{JiLa,JiLai}. See \cite{KLS} for details.

\gap
\it Remark. \rm The KLS theorem also provides an example indicating that 
the criterion of Theorem 1.1 is optimal in the sense that one cannot, in general,
say more by looking at the rate of growth of the $L^2$ norm. Indeed, by \eqref{5.4},
for a.e.\ $\omega,$ all solutions $\tilde u(n,E)$ for
every energy $E$ in the continuous spectrum satisfy 
\[ R^{-\rho}\|\tilde u (n,E)\|_{B_{R}}^{2} \leq C, \]
for any 
\[ \rho > 1+\frac{\lambda^2}{4-E^2} \]
and all $R.$ In particular, 
for every $\rho>1$ we can take $\lambda$ 
sufficiently small to ensure the existence of an interval $I_{\rho}$ around $E=0$
such that for a.e.\ $\omega$ all solutions (and in particular the one obeying
the boundary condition) satisfy
\[ R^{-\rho}\|\tilde u(n,E)\|_{B_{R}}^2 <C \]
for $E \in I_\rho$. Yet for a.e.\ $E \in I_\rho$, we have $D\mu(E)=0$ since the measure
is purely singular. This shows that no condition of type \eqref{1.2} with $\alpha >1$ 
leads in general to pointwise estimates on the derivatives of spectral measures.

This remark sounds trivial in one dimension, but it is straightforward --- using the analysis
of \cite{KLS} for the continuous analog $V_{\omega}(x)$ of 
the family of random potentials we study --- to give a similar example
which works in any dimension (in the continuous case). 
Set $H_{V_{\omega}} = -\Delta + \lambda V_{\omega}(r)$ with 
spherically symmetric potential. Using spherical symmetry, one shows 
that the spectrum of $H_{V_{\omega}}$ is purely singular with probability one. 
However, for every $\rho>1,$ there are solutions for a.e.\ $\omega$ 
and all energies $E$ sufficiently large such that \eqref{1.2} holds with 
$\alpha = \rho$. 

\gap
The following theorem shows that as long as our operators $h_\omega$ have
some continuous spectrum (which may be of arbitrarily small dimension),
their transport properties are arbitrarily close to ballistic.\\

\gap
\bf Theorem 5.1. \it Consider the family $h_{\omega}$ of random
Schr\"odinger operators defined on $\Z^{+}$
with potential $\lambda v_{\omega}(n)$,
where  $\lambda <2$ and the potential satisfies \eqref{5.3}.
Then for a.e.\ $\omega$,
for every $\psi$ such that $P_{c}(\omega) \psi \ne 0$
{\rm (}where $P_{c}(\omega)$ is the projector on the continuous spectrum
of the operator $h_{\omega}${\rm )} we have that for every $\epsilon > 0$
and $m > 0$
there is a positive constant $C_{\epsilon,m,\omega}$ such that for
any $T > 0$
\begin{equation}
\label{5.7}
\langle \langle |X|^{m} \psi(t), \psi(t) \rangle \rangle_{T}
\geq C_{\epsilon,m,\omega} T^{m(1-\epsilon)}.
\end{equation}

\gap
\begin{proof} \rm
By the results of the Gilbert-Pearson theory, the spectral measure
$\mu$ is supported on the set of the energies $E$ for which
the decaying solution \eqref{5.5} satisfies the boundary condition
(namely, $\theta(\omega)$ coincides with the Dirichlet boundary
condition). Moreover, these decaying solutions, which we will denote
by $u(n,E)$, are exactly the generalized eigenfunctions
in the sense of Theorem 1.2, if we normalize them by setting
$u(1,E)=1$.

Fix $\omega$ such that the results of the KLS theorem hold.
\eqref{5.5} implies that the generalized
eigenfunctions $u(n,E)$ of the operator $h_{\omega}$ satisfy
\begin{equation}
\label{5.9}
\limsup\limits_{R \rightarrow \infty} R^{-\gamma}\|u(n,E)\|_{B_{R}}^{2}
\leq \infty
\end{equation}
for every $\gamma > \alpha(E, \lambda)$ given by \eqref{5.4}.
Pick an open energy interval
$I = (E_1,E_2) \subset (-\sqrt{4 - \lambda^{2}}, \sqrt{4 -\lambda^{2}})$,
such that $0\notin I$, and $\mu^\psi(I) > 0$. Let
$\alpha_1 = \alpha(E_1, \lambda)$, $\alpha_2 = \alpha(E_2, \lambda)$.
$\alpha(E, \lambda)$ is monotone on $I$. Assume, without loss, that
$\alpha_1 < \alpha_2$. The restriction of $\mu^\psi$ to $I$ is
$\alpha_1$-continuous, and by \eqref{5.9},
$\limsup\limits_{R \rightarrow \infty}
R^{-\alpha_2}\|u(n,E)\|_{B_{R}}^{2} < \infty$ for any generalized
eigenfunction $u(n,E)$ with $E\in I$. Thus, by Theorem 1.2, for
each $m > 0$ there is a constant $C_{m,I,\omega}$ such that
for all $T > 0$
$$
\langle \langle |X|^{m} \psi(t), \psi(t) \rangle \rangle_{T}
\geq C_{m,I,\omega} T^{\frac{m\alpha_1}{\alpha_2}}.
$$
Since $P_c(\omega)\psi\not= 0$, we can clearly choose such an
interval $I$ with $\frac{\alpha_1}{\alpha_2} > 1-\epsilon$ and
$\mu^\psi(I) > 0$. Thus, Theorem 5.1 follows.
\end{proof}

\gap
\it Remark. \rm By using an extension of the proof of Theorem 1.2, one can
show that there is actually a component of the wave packet of size
corresponding to $\|P_c(\omega)\psi\|$ that is  spreading on average at
a rate which is arbitrarily close to ballistic. More explicitly, one
can show that for a.e.\ $\omega$, for every $\epsilon>0$ and $\rho >0$
there exists a constant $C_{\omega, \rho, \epsilon}$ such that if
$R_{T}= C_{\omega, \rho, \epsilon}T^{1-\epsilon}$, then
\begin{equation}
\label{5.8}
\langle \|P_{R_{T}}\psi (t) \|^{2} \rangle_{T} \leq
\|\psi - P_{c}(\omega)\psi \|^{2} + \rho.
\end{equation}
This easily yields Theorem 5.1 and is thus a stronger statement.

\gap
\begin{center} \bf Appendix 1. Generalizations \end{center}
\mgap
The whole proof of Theorem 1.1 readily extends to
more general settings. Namely, we can replace
 the operator $H_{V}^{\Omega}$ with general 
uniformly elliptic self-adjoint
operator $\Gamma$ such that
\[ \Gamma = (\partial_{l}-iA_{l}(x)) a_{lk}(x) (\partial_{k} - iA_{k}(x))+V(x) \]
provided that $a_{lk},$ $A_{l}$ and $V$ are ``nice enough'' (for example,
bounded and sufficiently smooth).
The proof for this case is very similar. The Green's formula leads
us to consider the following modified Wronskian:
\[  W_{\partial S}[f,g]  =  
 \int\limits_{\partial S} 
(\cos (\overline{n}, x_{l})a_{lk}((\partial_{k}-iA_{k})u)\overline{v} -
u \cos(\overline{n},x_{k})\overline{(\partial_{l} - iA_{l})a_{lk}v})d\sigma \]
It is clear that under our assumptions, the analog of Lemma 2.4
holds. The estimate of Lemma 2.3 also holds with the constant independent
of $R$ by the standard Sobolev estimates for bounded sufficiently
smooth coefficients (see, e.g., \cite{GT,Miz}).
The rest of the proof does not change.

A similar remark applies to some higher order operators and systems.
In particular, in one dimension, a self-adjoint half-line differential
operator of order $2n$ is given by the expression
\[ (Lf)(x) = (-1)^{n}(p_{0}f^{(n)})^{(n)}+ (-1)^{n-1}(p_{1}f^{(n-1)})^{(n-1)}+
\dots + p_{n}f \]
and a set of self-adjoint boundary conditions at zero.
The analog of the Wronskian in this case is determined by integration by
parts:
\begin{equation}
\label{A.1}
 W_{L,x}[f,g] = \sum\limits_{j=1}^{n} \sum\limits_{m=1}^{j}
(-1)^{m} \left( (p_{j}f^{(j)})^{(m-1)}g^{(j-m)}-f^{(j-m)}(p_{j}g^{(j)})^{(m-1)}
\right),
\end{equation}
where all values are taken at the point $x.$ The analog of the Sobolev
estimates of Lemma 2.3 is now the claim that for a solution of $(L-E)u =\phi,$
\[ \|u\|_{W_{2}^{m}(B_{R})} \leq C \|u\|_{L_{2}(B_{R+1})} \]
holds for
$m \leq 2n-1.$ Such estimates (in fact for $m \leq 2n$) are well-known
to hold for operators with bounded sufficiently smooth coefficients
(see, e.g., \cite{Miz}). The analog of Lemma 2.4 follows directly from
\eqref{A.1}; the rest of the proof of Theorem 1.1 does not change.

In particular, we have \\

\gap
\bf Theorem A.1. \it Let $L$ denote the self-adjoint differential operator
of order $2n$ with bounded sufficiently smooth
{\rm (}say, infinitely differentiable{\rm )} coefficients.
Suppose that for every $E$ in a set $S$ of positive
Lebesgue measure, there exists a bounded solution $u(x,E)$ of the generalized
eigenfunction equation
\[ (L -E)u =0 \]
satisfying the boundary conditions.
Suppose that for a compactly supported function $\phi \in L^{2},$ we
have
\[ \int u(x,E) \phi(x) dx \ne 0 \]
for a.e $E \in S.$
Then  the absolutely continuous part of the
spectral measure $\mu^{\phi}$
fills  $S$ {\rm(}so that $\mu^{\phi}(S_{1})>0$ for any $S_{1} \subset S$
of positive Lebesgue measure{\rm)}. \\

\gap
\it Remark. \rm Of course, we can also allow
for $\phi$ which are not
$L^{2}$, but from the Sobolev space $H_{-2}(H_{V}),$
such as the $\delta$ function and its
derivatives up to $2n-1,$ which are often used in the setting of
 one-dimensional differential
operators. The spectral measure is not finite in this case,
but nothing else changes. \\

\rm Theorem A.1 follows from the
above discussion and proof of Theorem 2.5.
This result may be viewed as a sort of an analog of \cite{Sim3,St}
for the higher order
case. It is typical, though, that our condition involves only one solution
(\cite{Sim3,St} requires all solutions to
be bounded) because the possible
multiplicity of the spectrum makes it unreasonable to demand all solutions
to be
bounded (in higher order cases) to get absolutely continuous spectrum.
On the other hand, our result does not guarantee pure absolute
continuity.

\gap
\begin{center} \bf Appendix 2. One-dimensional
perturbed Stark operators \end{center}
\mgap
In this appendix, we make a remark concerning dynamical
properties of a certain class of
perturbed Stark operators.
We denote by $H_{V,S}$ the operator defined on the whole
axis by the differential
expression
\[ -\frac{d^{2}}{dx^{2}} -x +V(x).  \]
Our results are based on the theorem proved in \cite{Kis}: \\

\gap
\bf Theorem. \it Suppose that $|V(x)|
\leq C(1+|x|)^{-\frac{1}{3}-\epsilon},$
or $V$ is bounded and has a derivative $V'$ which is bounded and H\"older continuous.
Then the whole axis $(-\infty, \infty)$ is an essential support of
the absolutely continuous part of the spectral measure $\mu.$ Moreover,
for a.e.\ $E \in R,$ there exist two linearly independent
solutions $u_{\pm}(x,E),$ such that
\[ u_{\pm}(x,E) = x^{-\frac{1}{4}}
\exp(\pm i(\frac{2}{3}x^{\frac{3}{2}}+f_{\pm}(x, E)))(1+o(1)) \]
as $x \rightarrow +\infty,$ where $|f_{\pm}'(x,E)| \leq C(1+x)^{-\frac{1}{2}}.$
\gap

\rm  Stark operators  do not fit into the framework
provided by Theorem 1.1 because of the strong
negative part of the potential (and resulting
failure of Lemma 2.3). Indeed, for a.e.\ energy $E$
here, we have a solution $u(x,E)$
which satisfies $R^{-\frac{1}{2}}\|u(x,E)\|_{B_{R}}\leq C(E),$ which, if
Theorem 1.1 were true,  would imply
$D^{\frac{1}{2}}(E) >0$ a.e.\ $E.$ It should be possible
 to prove an analog of Theorem 1.1 for
some perturbed Stark operators taking into account that
instead of the Sobolev estimates of
Lemma 2.3, one rather has
\[ \|\nabla u\|^{2}_{B_{R}} \leq CR\|u\|_{B_{R}}^{2}. \]

However, the criterion of Theorem 1.2 applies, giving immediately \\

\gap
\bf Theorem A.2. \it Under the conditions of the previous theorem,
for every vector $\psi$ with non-zero projection on
 the absolutely continuous subspace, we have
\[ \langle \langle |X|^{m}\psi(t), \psi(t) \rangle \rangle_{T} \geq CT^{2m}. \]

\gap
\rm  We note that there are examples
\cite{NP} of potentials $V$ satisfying
\[ |V(x)| \leq C(x)(1+|x|)^{-\frac{1}{3}}, \]
where $C(x)$ tends to infinity as $x \rightarrow \infty,$
but arbitrarily slowly,
such that for a corresponding Stark operator,
there is a dense set of eigenvalues embedded in the absolutely
continuous spectrum. Theorem A.2
shows that such potentials, nevertheless, do not
slow down dynamics corresponding to the absolutely continuous component. \\

\gap
\bf Acknowledgements. 
\rm We thank Y.~Avron, I.~Guarneri, R.~Ketzmerick, B.~Simon and S.~Tcheremchantsev
for stimulating
discussions. We are grateful to the referees for useful sugestions and 
corrections.
AK's research is supported in part by NSF grant DMS-9801530. YL's
research is supported in part by NSF grant DMS-9801474.\\

\gap\\
\noindent Department of Mathematics \\
University of Chicago \\
5734 S. University Ave. \\
Chicago, IL 60637 \\
email: kiselev@math.uchicago.edu \\

\noindent and \\

\noindent Department of Mathematics  \\
California Institute of Technology \\
Pasadena, CA 91125 \\
email: ylast@cco.caltech.edu

\end{document}